\documentclass[11pt]{amsart}
\usepackage[margin=2.7cm]{geometry}
\usepackage{amsmath,amssymb,amsthm,booktabs,enumitem}
\usepackage[hidelinks]{hyperref}

\newtheorem{thm}{Theorem}[section]
\newtheorem{prop}[thm]{Proposition}
\newtheorem{lem}[thm]{Lemma}
\newtheorem{cor}[thm]{Corollary}
\theoremstyle{definition}

\newtheorem{ex}[thm]{Example}
\theoremstyle{remark}
\newtheorem{rk}[thm]{Remark}

\newcommand{\C}{\mathbb{C}}
\newcommand{\PP}{\mathbb{P}}
\newcommand{\Z}{\mathbb{Z}}
\newcommand{\AR}{\mathrm{AR}}
\newcommand{\Sp}{\mathrm{Sp}}
\newcommand{\Gr}{\mathrm{Gr}}
\newcommand{\cO}{\mathcal{O}}
\newcommand{\cE}{\mathcal{E}}
\DeclareMathOperator{\rank}{rank}

\begin{document}

\title[Symmetry of pole order spectra and freeness]{Symmetry of pole order spectra and freeness of plane curves}

\author{Gabriel Sticlaru}
\address{Faculty of Mathematics and Informatics, Ovidius University, Bd. Mamaia 124, 900527 Constan\c ta, Romania}
\email{gabriel.sticlaru@gmail.com}


\begin{abstract}
Let $C:f=0$ be a reduced plane curve of degree $d$ having only weighted homogeneous singularities, and let $\Sp_P^0(f)$ be the pole order spectrum of the top cohomology of its Milnor fiber. Writing $\Sp^0_P(f)=\sum_q m_q\,t^{q/d}$, we show that
$m_q-m_{2d-q}=n(f)_{q-3}-n(f)_{q+d-3}$ for every integer $q$, where $n(f)_k$ is the Hilbert function of the Jacobian module $N(f)$. Consequently, $\Sp^0_P(f)$ is symmetric with respect to $\alpha\mapsto 2-\alpha$ if and only if $C$ is free. For a free curve with exponents $(d_1,d_2)$ we give a closed formula: $\Sp^0_P(f)-t^{-1}\Sp^1_P(f)$ depends only on $(d,d_1)$. This explains the symmetry of pole order spectra of reflection line arrangements, an open problem raised in Dimca's monograph on hyperplane arrangements.
\end{abstract}

\maketitle

\section{Introduction}

Let $S=\C[x,y,z]$ and let $C:f=0$ be a reduced curve of degree $d\ge 2$ in $\PP^2$. The Milnor fiber $F:f(x,y,z)=1$ carries the monodromy action and, on its cohomology, the pole order filtration $P$. Dimca and the author gave algorithms computing the corresponding pole order spectra $\Sp^0_P(f)$ and $\Sp^1_P(f)$ (of $H^2(F,\C)$ and $H^1(F,\C)$) via a spectral sequence whose $E_1$-term is built from the Koszul complex of the partial derivatives of $f$ \cite{DS15,DS16a,DS16,DS17,Di17}.

In all the examples coming from reflection groups, the spectrum $\Sp^0_P(f)$ turned out to be symmetric: the coefficient of $t^\alpha$ equals that of $t^{2-\alpha}$. The prototype is the arrangement $\mathcal A(2,2,3):(x^2-y^2)(y^2-z^2)(z^2-x^2)=0$, for which \cite[Example~8.14]{Di17}
\begin{equation}\label{eq:A223}
\Sp^0_P(f)=t^{3/6}+3t^{4/6}+2t^{5/6}+6t+2t^{7/6}+3t^{8/6}+t^{9/6}.
\end{equation}
This symmetry is ``very mysterious'', since the corresponding Hodge spectra are not symmetric, and explaining it is stated as an open problem in \cite[Remark~8.13]{Di17}; see also \cite[Remark~5.6]{DS17}. The aim of this note is to settle the question for plane curves.

For $k\in\Z$ let $n(f)_k=\dim N(f)_k$, where $N(f)=\widehat J_f/J_f$ is the Jacobian module \cite{Se}, $J_f=(f_x,f_y,f_z)$ and $\widehat J_f$ its saturation with respect to $\mathfrak m=(x,y,z)$. All exponents of $\Sp^0_P(f)$ lie in $\frac1d\Z_{>0}$, so we write
\[
\Sp_P^0(f)=\sum_{q\ge1} m_q\, t^{q/d} .
\]

\begin{thm}\label{thmA}
Let $C:f=0$ be a reduced plane curve of degree $d$ with only weighted homogeneous singularities. Then for every $q\in\Z$
\begin{equation}\label{eq:defect}
m_q-m_{2d-q}=n(f)_{q-3}-n(f)_{q+d-3}.
\end{equation}
In particular, $\Sp^0_P(f)$ is symmetric with respect to $\alpha\mapsto 2-\alpha$ if and only if $C$ is a free curve.
\end{thm}

All line arrangements have only weighted homogeneous (ordinary) singularities, and all reflection arrangements are free by Terao's theorem \cite{OT}. Hence:

\begin{cor}\label{corB}
The pole order spectrum $\Sp^0_P(f)$ of any reflection line arrangement in $\PP^2$ (for instance $\mathcal A(m,m,3)$ and $\mathcal A(m,1,3)$, $m\ge 2$, or the arrangements $A_3$, $B_3$, $H_3$) is symmetric. More generally, a line arrangement has a symmetric spectrum $\Sp^0_P(f)$ if and only if it is free.
\end{cor}

Since freeness of line arrangements is the subject of Terao's conjecture, Corollary~\ref{corB} has the following reformulation: \emph{Terao's conjecture holds for line arrangements in $\PP^2$ if and only if the symmetry of $\Sp^0_P(f)$ is determined by the intersection lattice.} The spectrum $\Sp^0_P(f)$ itself is not combinatorial \cite[Remark~8.12]{Di17}.

For free curves we obtain an explicit formula. For $n\in\Z$ set $s(n)=\binom{n+2}{2}$ if $n\ge0$ and $s(n)=0$ otherwise, so that $s(n)=\dim S_n$.

\begin{thm}\label{thmC}
Let $C:f=0$ be a free curve of degree $d$ with exponents $d_1\le d_2$, $d_1+d_2=d-1$, having only weighted homogeneous singularities. Then
\[
\Sp^0_P(f)=\sum_{q} G_{d,d_1}(q)\,t^{q/d}\;+\;t^{-1}\,\Sp^1_P(f),
\]
where the integers $G_{d,d_1}(q)$ are given by the generating polynomial
\begin{equation}\label{eq:Gfree}
\sum_q G_{d,d_1}(q)T^q=\frac{T^3-T^{d_1+2}-T^{d_2+2}+3T^{d+1}-3T^{d+2}+T^{d+d_1+2}+T^{d+d_2+2}-T^{2d}}{(1-T)^3}.
\end{equation}
This is a palindromic polynomial of degree $2d-3$, with $T^{2d}G_{d,d_1}(1/T)=G_{d,d_1}(T)$ and $G_{d,d_1}(1)=d(d_1-1)(d_2-1)-1$.
\end{thm}

In other words, the ``modified spectrum'' $\Sp^0_P(f)-t^{-1}\Sp^1_P(f)$ of a free curve depends only on $d$ and $d_1$. This is the plane curve analogue of a symmetry proved by M.~Saito for a modified pole order spectrum of strongly free divisors in $\PP^3$ \cite[Theorem~4.5]{Sa17}. The new points here are that, for plane curves, the term $t^{-1}\Sp^1_P(f)$ is itself symmetric (Proposition~\ref{prop:eps}), so the symmetry passes to the true spectrum $\Sp^0_P(f)$, and that the failure of symmetry for non-free curves is measured exactly by the Jacobian module.

Saito expressed the pole order spectrum of a projective hypersurface with weighted homogeneous isolated singularities through the Hilbert series of the $E_2$-terms, and proved that the torsion part of the Milnor algebra, which in our case is $N(f)$, injects into $E_2$ \cite[Proposition~1 and Theorem~3]{Sa16}. He did not consider the symmetry $\alpha\mapsto 2-\alpha$. Our contribution is to identify the defect of this symmetry with $N(f)$, using the symmetry of the $H^1$-contribution and Serre duality for the syzygy bundle.

The paper is organized as follows. Section~2 recalls the spectral sequence. Section~3 proves the key identity $m_q=G_f(q)+\varepsilon_q$. Section~4 shows that the $H^1$-contribution $\varepsilon_q$ is symmetric. Section~5 proves Theorems~\ref{thmA} and \ref{thmC} using Serre duality for the syzygy bundle. Section~6 contains examples.

\section{The pole order spectral sequence}

Let $\Omega^j$ be the graded $S$-module of global polynomial $j$-forms on $\C^3$, with $\deg x=\deg dx=1$, so that $\Omega^j_q\cong S_{q-j}^{\binom 3j}$. Let $K^*_f=(\Omega^*,df\wedge)$ be the Koszul complex of $f_x,f_y,f_z$. Since $C$ is reduced, its singularities are isolated and
\begin{equation}\label{eq:koszul}
H^0(K^*_f)=H^1(K^*_f)=0,\qquad H^3(K^*_f)_q\cong M(f)_{q-3},
\end{equation}
where $M(f)=S/J_f$ is the Milnor algebra. Moreover $H^2(K^*_f)_q=\mathrm{Syz}(f)_q/(df\wedge\Omega^1)_q$, where $\mathrm{Syz}(f)_q=\ker\{df\wedge:\Omega^2_q\to\Omega^3_{q+d}\}\cong \AR(f)_{q-2}$ and $\AR(f)$ is the graded $S$-module of Jacobian syzygies $(a,b,c)$, $af_x+bf_y+cf_z=0$ \cite[Section~8.3]{Di17}.

\begin{thm}[\cite{DS16,Di17}]\label{thm:ss}
For each integer $1\le k\le d$ there is a spectral sequence $E_*(f)_k$ with
\[
E_1^{s,t}(f)_k=H^{s+t+1}(K^*_f)_{td+k}\ \Longrightarrow\ E_\infty^{s,t}(f)_k=\Gr_P^s\widetilde H^{s+t}(F,\C)_\lambda,\qquad \lambda=\exp(-2\pi i k/d).
\]
The differential $d_1$ is induced by the exterior derivative, $d':H^2(K^*_f)\to H^3(K^*_f)$.
\end{thm}

We use the following known facts, valid when $C$ has only weighted homogeneous singularities.

\begin{enumerate}[label=(F\arabic*)]
\item\label{F1} (M.~Saito) $E_2(f)_k=E_\infty(f)_k$ for all $k$ \cite{Sa16}, see also \cite[Theorem~8.13]{Di17}. For $k=d$, i.e.\ for the complement $U$, this was proved earlier in \cite[Theorem~2.4]{DS15}.
\item\label{F2} $\dim\ker\{d':H^2(K^*_f)_q\to H^3(K^*_f)_q\}$ is bounded by the number of spectral numbers of the singularities of $C$ equal to $q/d$ \cite[Theorem~8.14]{Di17}. Since these numbers lie in $(0,2)$, the kernel vanishes for $q\ge 2d$.
\item\label{F3} For $1\le k\le d-1$ the multiplicity of $\lambda=\exp(-2\pi ik/d)$ as a root of the Alexander polynomial, i.e. $\dim H^1(F,\C)_\lambda$, equals $\dim E_2^{1,0}(f)_k+\dim E_2^{1,0}(f)_{d-k}$ \cite[Theorem~8.15]{Di17}.
\end{enumerate}

For $q\in\Z$ write uniquely $q=td+k$ with $1\le k\le d$, and set
\[
m'_q=\dim E_2^{2-t,t}(f)_k,\qquad \varepsilon_q=\dim E_2^{1-t,t}(f)_k,
\]
with the convention that these are $0$ whenever the corresponding term is not in the first quadrant. By Theorem~\ref{thm:ss} and \ref{F1}, $m'_q=\dim\Gr_P^{2-t}H^2(F,\C)_\lambda$ and $\varepsilon_q=\dim \Gr_P^{1-t}H^1(F,\C)_\lambda$ (for $\lambda\neq1$; for $\lambda=1$ we use reduced cohomology). Recall \cite[(8.44)]{Di17} that $\Sp^j_P(f)=\sum_{\alpha>0}m^j_\alpha t^\alpha$ with $m^j_\alpha=\dim\Gr^p_PH^{2-j}(F,\C)_{\exp(-2\pi i\alpha)}$, $p=[3-\alpha]$. A direct check of the indices shows that this agrees with the notation of the Introduction:

\begin{lem}\label{lem:index}
$m'_q=m_q$ is the coefficient of $t^{q/d}$ in $\Sp^0_P(f)$, and $\varepsilon_q$ is the coefficient of $t^{(q+d)/d}$ in $\Sp^1_P(f)$. In other words
\[
\Sp^0_P(f)=\sum_q m_q\,t^{q/d},\qquad t^{-1}\Sp^1_P(f)=\sum_q \varepsilon_q\, t^{q/d}.
\]
\end{lem}

\begin{proof}
If $q=td+k$, then $\alpha=q/d=t+k/d$ satisfies $\exp(-2\pi i\alpha)=\lambda$ and $[3-\alpha]=2-t$, which is the filtration index of $E_2^{2-t,t}$. Similarly $\alpha=(q+d)/d=t+1+k/d$ gives $[3-\alpha]=1-t$, the filtration index of $E^{1-t,t}_2$ in $H^1$.
\end{proof}

\section{The key identity}

Let $\AR_f(T)=\sum_j\dim\AR(f)_jT^j$ be the Hilbert series of $\AR(f)$ and set
\[
A_d(T)=\frac{T^3-3T^{d+2}+3T^{d+1}-T^{2d}}{(1-T)^3}.
\]

\begin{prop}\label{prop:key}
For every $q\in\Z$ one has $m_q=G_f(q)+\varepsilon_q$, where
\begin{equation}\label{eq:Gf}
G_f(q)=\dim M(f)_{q-3}-\dim\AR(f)_{q-2}+\dim(df\wedge\Omega^1)_q .
\end{equation}
Moreover $\sum_qG_f(q)T^q=A_d(T)-T^2(1-T^d)\,\AR_f(T)$.
\end{prop}

\begin{proof}
In the spectral sequence $E_*(f)_k$, the term $E_1^{2-t,t}=H^3(K^*_f)_q$ has outgoing differential landing in $H^4=0$ and incoming differential $d_1=d':H^2(K^*_f)_q\to H^3(K^*_f)_q$, whose kernel is $E_2^{1-t,t}$. Hence
\[
m_q=\dim H^3(K^*_f)_q-\dim H^2(K^*_f)_q+\varepsilon_q .
\]
Now \eqref{eq:koszul} gives $\dim H^3(K^*_f)_q=\dim M(f)_{q-3}$ and $\dim H^2(K^*_f)_q=\dim\AR(f)_{q-2}-\dim(df\wedge\Omega^1)_q$, which proves \eqref{eq:Gf}.

For the generating function, the exactness of $K^*_f$ in degrees $0$ and $1$ yields the exact sequence $0\to\Omega^0_{q-2d}\to\Omega^1_{q-d}\to(df\wedge\Omega^1)_q\to0$, hence $\dim(df\wedge\Omega^1)_q=3s(q-d-1)-s(q-2d)$. Since $J_f$ is the image of $S(1-d)^3\to S$ with kernel $\AR(f)(1-d)$, we get $\dim M(f)_j=s(j)-3s(j-d+1)+\dim\AR(f)_{j-d+1}$. Substituting into \eqref{eq:Gf} and using $\sum_q s(q-a)T^q=T^a/(1-T)^3$ gives the formula.
\end{proof}

\begin{lem}\label{lem:Asym}
$T^{2d}A_d(1/T)=A_d(T)$.
\end{lem}

\begin{proof}
Since $(1-T^{-1})^3=-T^{-3}(1-T)^3$, the claim is equivalent to $-T^{2d+3}N(1/T)=N(T)$ for $N(T)=T^3-3T^{d+2}+3T^{d+1}-T^{2d}$, which is immediate.
\end{proof}

\section{Symmetry of the $H^1$-contribution}

\begin{prop}\label{prop:eps}
If $C$ has only weighted homogeneous singularities, then $\varepsilon_q=\varepsilon_{2d-q}$ for all $q\in\Z$.
\end{prop}

\begin{proof}
By definition $\varepsilon_q=0$ for $q\le0$, and $\varepsilon_q=0$ for $q\ge 2d$ by \ref{F2} (for $q>2d$ also because $E_2(f)_k$ lies in the first quadrant \cite[Corollary~8.3]{Di17}). So it remains to treat $1\le q\le 2d-1$, and by the symmetry of the statement we may assume $d+1\le q\le 2d-1$, i.e. $q=d+k$ with $1\le k\le d-1$. Since $E_2=E_\infty$ by \ref{F1}, the abutment of $E_*(f)_k$ in total degree $1$, together with the fact that $E_2(f)_k$ lies in the first quadrant, gives
\[
\dim H^1(F,\C)_\lambda=\dim E_2^{1,0}(f)_k+\dim E_2^{0,1}(f)_k=\varepsilon_k+\varepsilon_{d+k},
\]
where $\lambda=\exp(-2\pi i k/d)$. On the other hand, \ref{F3} gives $\dim H^1(F,\C)_\lambda=\varepsilon_k+\varepsilon_{d-k}$. Hence $\varepsilon_{d+k}=\varepsilon_{d-k}$, which is the claim for $q=d+k$.
\end{proof}

\begin{rk}
For $q=d$ (i.e. $\lambda=1$) the statement is empty. Since $\varepsilon_{2d}=0$, we have $P^1H^1(F,\C)_1=H^1(F,\C)_1$ and $\varepsilon_d=\dim H^1(F,\C)_1=\dim H^1(U)=r-1$, where $r$ is the number of irreducible components of $C$ and $U=\PP^2\setminus C$; compare \cite[Remark~2.3 and Section~3]{DS15}.
\end{rk}

\section{Proofs of the main results}

\subsection{The syzygy bundle}
Let $\cE$ be the kernel of the sheaf morphism $\cO^3\to\cO(d-1)$, $(a,b,c)\mapsto af_x+bf_y+cf_z$. Its image is $\mathcal I_Z(d-1)$, where $Z$ is the zero-dimensional Jacobian subscheme of $\PP^2$, so $\cE$ is a reflexive sheaf of rank $2$ on a smooth surface, i.e. a vector bundle, with $c_1(\cE)=1-d$. Thus $\cE^\vee\cong\cE(d-1)$. Since $\AR(f)$ is a second syzygy module, it has depth $\ge 2$, hence
\[
H^0(\cE(j))=\AR(f)_j\quad\text{for all } j.
\]
Taking global sections in $0\to\cE(j)\to\cO(j)^3\to\mathcal I_Z(j+d-1)\to0$ and using $H^0(\mathcal I_Z(k))=(\widehat J_f)_k$ and $H^1(\cO(j))=0$ we get
\begin{equation}\label{eq:H1}
H^1(\cE(j))\cong (\widehat J_f/J_f)_{j+d-1}=N(f)_{j+d-1}.
\end{equation}
Finally, Serre duality gives $H^2(\cE(j))\cong H^0(\cE^\vee(-j-3))^\vee=H^0(\cE(d-4-j))^\vee$.

\begin{lem}\label{lem:duality}
As rational functions,
$\AR_f(T)+T^{d-4}\AR_f(1/T)=\sum_j n(f)_{j+d-1}T^j.$
\end{lem}

\begin{proof}
Let $\chi(j)=\chi(\cE(j))$, a polynomial in $j$. For every $j$ we have $h^0(\cE(j))+h^2(\cE(j))=\chi(j)+h^1(\cE(j))$, and $\sum_jh^2(\cE(j))T^j=\sum_jh^0(\cE(d-4-j))T^j=T^{d-4}\AR_f(1/T)$. Choose $J\gg0$ such that $h^1(\cE(j))=0$ for $|j|\ge J$, $h^2(\cE(j))=0$ for $j\ge J$, and $h^0(\cE(j))=0$ for $j\le -J$. Then
\[
\AR_f(T)+T^{d-4}\AR_f(1/T)-\sum_jh^1(\cE(j))T^j=\sum_{j\ge J}\chi(j)T^j+\sum_{j\le -J}\chi(j)T^j+\sum_{|j|<J}\chi(j)T^j .
\]
For a polynomial $\chi$, the rational functions represented by $\sum_{j\ge J}\chi(j)T^j$ (expanded at $0$) and $\sum_{j<J}\chi(j)T^j$ (expanded at $\infty$) add up to zero. Hence the right-hand side vanishes, and \eqref{eq:H1} concludes the proof.
\end{proof}

\begin{rk}\label{rk:selfdual}
Replacing $T$ by $1/T$ and multiplying by $T^{d-4}$ leaves the left-hand side of Lemma~\ref{lem:duality} unchanged, hence also the right-hand side. This gives $n(f)_k=n(f)_{3(d-2)-k}$, the known self-duality of the Jacobian module \cite{Se}, see also \cite{DSa}. So Lemma~\ref{lem:duality} contains this symmetry, and it is what makes the defect in Theorem~\ref{thmA} antisymmetric under $q\mapsto 2d-q$, as it must be.
\end{rk}

\subsection{Proof of Theorem~\ref{thmA}}
Since $m_q$ and $\varepsilon_q$ vanish for $q\le 0$ and $q>3d$, $G_f(T)=\sum_qG_f(q)T^q$ is a Laurent polynomial. By Proposition~\ref{prop:key}, Lemma~\ref{lem:Asym} and $T^{2d}\cdot T^{-2}(1-T^{-d})=-T^{d-2}(1-T^d)$,
\[
G_f(T)-T^{2d}G_f(1/T)=-T^2(1-T^d)\bigl(\AR_f(T)+T^{d-4}\AR_f(1/T)\bigr)=-T^2(1-T^d)\sum_jn(f)_{j+d-1}T^j,
\]
using Lemma~\ref{lem:duality}. Comparing coefficients of $T^q$ gives $G_f(q)-G_f(2d-q)=n(f)_{q-3}-n(f)_{q+d-3}$. By Proposition~\ref{prop:eps}, $m_q-m_{2d-q}=G_f(q)-G_f(2d-q)$, and Lemma~\ref{lem:index} gives \eqref{eq:defect}.

If $C$ is free, then $\AR(f)$ is free, so $M(f)$ has projective dimension $2$ and, by the Auslander--Buchsbaum formula, depth $1$; hence $N(f)=H^0_{\mathfrak m}(M(f))=0$ and the spectrum is symmetric (that $C$ is free if and only if $N(f)=0$ is well known, see \cite{Se,DS-RIMS}). Conversely, if $\Sp^0_P(f)$ is symmetric, then $n(f)_{k}=n(f)_{k+d}$ for all $k$. Since $N(f)$ is a module of finite length, $n(f)_k=0$ for $|k|\gg0$, and periodicity forces $N(f)=0$. Then $M(f)$ has depth $1$, $\AR(f)$ has projective dimension $0$, and $C$ is free. \qed

\subsection{Proof of Theorem~\ref{thmC}}
If $C$ is free with exponents $(d_1,d_2)$, then $\AR(f)\cong S(-d_1)\oplus S(-d_2)$, so $\AR_f(T)=(T^{d_1}+T^{d_2})/(1-T)^3$. Proposition~\ref{prop:key} gives $G_f=G_{d,d_1}$ as in \eqref{eq:Gfree}, and Lemma~\ref{lem:index} gives the formula for $\Sp^0_P(f)$. The palindromic property follows from Theorem~\ref{thmA} with $N(f)=0$, or directly: under $T\mapsto 1/T$ the numerator of \eqref{eq:Gfree}, multiplied by $-T^{2d+3}$, is permuted into itself, using $d_1+d_2=d-1$.

Since $G_f$ is a Laurent polynomial starting with $T^3$, the palindromic property shows that $G_{d,d_1}$ is a polynomial supported in $[3,2d-3]$. Finally, $G_{d,d_1}(1)=\sum_qm_q-\sum_q\varepsilon_q=\dim H^2(F)-\dim H^1(F)=\chi(F)-1$ (reduced cohomology), and $\chi(F)=d\,\chi(U)$. Now $\chi(U)=3-\chi(C)=3-(3d-d^2+\mu(C))$, and for weighted homogeneous singularities $\mu(C)=\tau(C)$, while $\tau(C)=(d-1)^2-d_1d_2$ for a free curve. Using $d_1+d_2=d-1$ we get $\chi(U)=(d_1-1)(d_2-1)$, hence $G_{d,d_1}(1)=d(d_1-1)(d_2-1)-1$. \qed

\begin{rk}
The weighted homogeneity assumption cannot simply be dropped. According to the computations in \cite{DS16}, an irreducible free sextic with exponents $(2,3)$ and a singularity which is not weighted homogeneous has
$\Sp^0_P(f)=t^{3/6}+t^{4/6}+2t^{5/6}+2t^{7/6}+t^{8/6}$, which is not symmetric. For such curves the facts \ref{F1}--\ref{F3} are not available, and the identities above, which concern $E_2$-terms, need not compute the spectrum.
\end{rk}

\section{Examples and computations}

All the spectra below were computed by implementing the formula $m_q=\dim M(f)_{q-3}-\rank\delta_q$, where $\delta_q:\mathrm{Syz}(f)_q\to M(f)_{q-3}$ sends $(a,b,c)$ to the class of $a_x+b_y+c_z$ \cite[(8.36)--(8.37)]{Di17}. The computations were carried out independently with a Python implementation (linear algebra modulo $32003$) and with Singular \cite{Sing}; the two give identical results. The Singular scripts are available from the author upon request.

\begin{ex}[Reflection arrangements]
For $\mathcal A(2,2,3)$ ($d=6$, exponents $(2,3)$) one gets $G_{6,2}=T^3+2T^4+2T^5+T^6+2T^7+2T^8+T^9$ and $t^{-1}\Sp^1_P(f)=t^{4/6}+5t+t^{8/6}$. Their sum is exactly \eqref{eq:A223}. Table~\ref{tab:refl} lists the data for the monomial arrangements $\mathcal A(m,m,3)$ and $\mathcal A(m,1,3)$ for small $m$. In all cases $\Sp^0_P(f)$ is symmetric, $m_q-G_{d,d_1}(q)=\varepsilon_q$ is supported in $\{2d/3,\,d,\,4d/3\}$ (cubic roots of unity, compare \cite[Theorem~8.18]{Di17}), and $\varepsilon_d=d-1$.
\end{ex}

\begin{table}[ht]
\centering
\begin{tabular}{lcccl}
\toprule
arrangement & $d$ & $(d_1,d_2)$ & $\dim H^2(F)$ & nonzero $\varepsilon_q$ \\
\midrule
$\mathcal A(2,2,3)$ & 6 & (2,3) & 18 & $\varepsilon_4=1,\ \varepsilon_6=5,\ \varepsilon_8=1$\\
$\mathcal A(3,3,3)$ & 9 & (4,4) & 92 & $\varepsilon_6=2,\ \varepsilon_9=8,\ \varepsilon_{12}=2$\\
$\mathcal A(4,4,3)$ & 12 & (5,6) & 252 & $\varepsilon_8=1,\ \varepsilon_{12}=11,\ \varepsilon_{16}=1$\\
$\mathcal A(5,5,3)$ & 15 & (6,8) & 540 & $\varepsilon_{10}=1,\ \varepsilon_{15}=14,\ \varepsilon_{20}=1$\\
$\mathcal A(2,1,3)=B_3$ & 9 & (3,5) & 79 & $\varepsilon_9=8$\\
$\mathcal A(3,1,3)$ & 12 & (4,7) & 226 & $\varepsilon_{12}=11$\\
$\mathcal A(4,1,3)$ & 15 & (5,9) & 495 & $\varepsilon_{10}=1,\ \varepsilon_{15}=14,\ \varepsilon_{20}=1$\\
\bottomrule
\end{tabular}
\caption{Free reflection line arrangements: $\Sp^0_P(f)=G_{d,d_1}+\sum_q\varepsilon_q t^{q/d}$, all symmetric.}\label{tab:refl}
\end{table}

\begin{ex}[$\mathcal A(3,3,3)$]
Here $d=9$, $(d_1,d_2)=(4,4)$, and
\[
\Sp^0_P(f)=t^{3/9}+3t^{4/9}+6t^{5/9}+10t^{6/9}+9t^{7/9}+9t^{8/9}+16t+9t^{10/9}+9t^{11/9}+10t^{12/9}+6t^{13/9}+3t^{14/9}+t^{15/9},
\]
while $G_{9,4}(T)=T^3+3T^4+6T^5+8T^6+9T^7+9T^8+8T^9+9T^{10}+9T^{11}+8T^{12}+6T^{13}+3T^{14}+T^{15}$. The difference $2t^{6/9}+8t+2t^{12/9}=t^{-1}\Sp^1_P(f)$ reflects the two-dimensional eigenspaces $H^1(F,\C)_\lambda$ for $\lambda$ a primitive cubic root of unity.
\end{ex}

\begin{ex}[Free arrangements with trivial monodromy]
If $\mathcal A$ is a free arrangement of $d$ lines with $H^1(F,\C)_{\neq1}=0$, then $\Sp^1_P(f)=(d-1)t^2$ and Theorem~\ref{thmC} gives $\Sp^0_P(f)=G_{d,d_1}(t^{1/d})+(d-1)t$. Thus $\Sp^0_P(f)$ depends only on $d$ and $d_1$, and in particular it is determined by the intersection lattice. For example, the free arrangement $xyz(x-y)(x-z)(y-z)(x+y-z)$ with exponents $(3,3)$ has
$\Sp^0_P(f)=t^{3/7}+3t^{4/7}+4t^{5/7}+4t^{6/7}+9t+4t^{8/7}+4t^{9/7}+3t^{10/7}+t^{11/7}$.
This is in contrast with the two arrangements of Ziegler's pair: they have the same intersection lattice, are not free, and have different spectra $\Sp^0_P$ \cite[Remark~8.12]{Di17}. By Theorem~\ref{thmA}, the difference is detected by their Jacobian modules.
\end{ex}

\begin{ex}[Non-free curves]
For non-free curves, Theorem~\ref{thmA} measures the lack of symmetry. Table~\ref{tab:nonfree} compares the observed defect with the Hilbert function of $N(f)$. In each case \eqref{eq:defect} was verified for all $q$.
\end{ex}

\begin{table}[ht]
\centering
\begin{tabular}{llll}
\toprule
curve & $d$ & $\Sp^0_P(f)$ & nonzero $n(f)_k$ \\
\midrule
5 general lines & 5 & $t^{3/5}+3t^{4/5}+6t+3t^{6/5}+3t^{7/5}+2t^{8/5}$ & $n_4=n_5=2$\\
$xyz(x+y+z)(x-y)$ & 5 & $t^{3/5}+2t^{4/5}+5t+2t^{6/5}+2t^{7/5}+t^{8/5}$ & $n_4=n_5=1$\\
$y^4z-x^5$ (nearly free) & 5 & $t^{6/5}+t^{7/5}+t^{8/5}+t^{9/5}$ & $n_3=\dots=n_6=1$\\
6 general lines & 6 & $t^{3/6}+3t^{4/6}+6t^{5/6}+10t+6t^{7/6}+6t^{8/6}+5t^{9/6}+3t^{10/6}$ & $n_5=3,\ n_6=4,\ n_7=3$\\
\bottomrule
\end{tabular}
\caption{Non-free curves: the asymmetry of $\Sp^0_P(f)$ equals $n(f)_{q-3}-n(f)_{q+d-3}$.}\label{tab:nonfree}
\end{table}

\begin{rk}[Higher dimensions]
Corollary~\ref{corB} explains the symmetry for reflection arrangements of rank $3$. The examples in \cite{DS17} concern rank-$4$ reflection arrangements (surfaces in $\PP^3$), where some free, decomposable arrangements $(x^p+y^p)(z^q+w^q)=0$ have non-symmetric spectra. There the spectral sequence is only known to degenerate at $E_3$ \cite{Sa19}, and Saito's symmetry \cite[Theorem~4.5]{Sa17} concerns a modified spectrum. The method above suggests the following question: for an essential free arrangement in $\PP^3$, is the $H^2$-contribution to Saito's modified spectrum symmetric, as $\varepsilon_q$ is for plane curves?
\end{rk}

\section*{Acknowledgments}
The author is grateful to Alexandru Dimca for reading this paper and for many very useful suggestions.

\end{document}